\documentclass[10pt,twocolumn,twoside]{IEEEtran}

\usepackage{cite}
\usepackage{amsmath,amssymb,amsfonts}
\newtheorem{remark}{Remark}

\usepackage{algorithmic}
\usepackage{graphicx}
\usepackage{textcomp}
\usepackage{subcaption}

\def\BibTeX{{\rm B\kern-.05em{\sc i\kern-.025em b}\kern-.08em
    T\kern-.1667em\lower.7ex\hbox{E}\kern-.125emX}}

\begin{document}

\title{A Deep Reinforcement Learning Approach to Efficient Distributed Optimization}
\author{Daokuan Zhu, Tianqi Xu, and Jie Lu, \IEEEmembership{Member, IEEE}
\thanks{D. Zhu is with the School of Information Science and Technology, ShanghaiTech University, Shanghai 201210, China, with the University of Chinese Academy of Sciences, Beijing 100049, China, and with the Shanghai Institute of Microsystem and Information Technology, Chinese Academy of Sciences, Shanghai 200050, China. Email: \texttt{\small zhudk@shanghaitech.edu.cn}.}
\thanks{T. Xu and J. Lu are with the School of Information Science and Technology, ShanghaiTech University, Shanghai 201210, China. Email: \texttt{\small \{xutq2022, lujie\}@shanghaitech.edu.cn}.}}

\maketitle

\begin{abstract}
In distributed optimization, the practical problem-solving performance is essentially sensitive to algorithm selection, parameter setting, problem type and data pattern.
Thus, it is often laborious to acquire a highly efficient method for a given specific problem.
In this paper, we propose a learning-based method to achieve efficient distributed optimization over networked systems.
Specifically, a deep reinforcement learning (DRL) framework is developed for adaptive configuration within a parameterized unifying algorithmic form, which incorporates an abundance of decentralized first-order and second-order optimization algorithms.
We exploit the local consensus and objective information to represent the regularities of problem instances and trace the solving progress, which constitute the states observed by a DRL agent.
The framework is trained using Proximal Policy Optimization (PPO) on a number of practical problem instances of similar structures yet different problem data.
Experiments on various smooth and non-smooth classes of objective functions demonstrate that our proposed learning-based method outperforms several state-of-the-art distributed optimization algorithms in terms of convergence speed and solution accuracy.
\end{abstract}


\newcommand{\xb}{\mathbf{x}}
\newcommand{\Rbb}{\mathbb{R}}
\newcommand{\qb}{\mathbf{q}}
\newcommand{\vb}{\mathbf{v}}

\newcommand{\st}{\text{s.t.}}
\newcommand{\argminx}{\underset{x \in \mathbb{R}^d}{\arg\min}\ }
\newcommand{\argminxb}{\underset{\mathbf{x} \in \mathbb{R}^{Nd}}{\arg\min}\ }

\newcommand{\sumi}{\sum_{i=1}^N}

\newcommand{\alphaik}{\alpha_i^k}
\newcommand{\betaik}{\beta_i^k}
\newcommand{\rhoik}{\rho_i^k}

\section{Introduction}
\label{sec:intro}
With rapid developments of networked systems such as Internet of Vehicles (IoV), smart grids, smart buildings and fog computing, applying traditional centralized optimization methods in these systems faces the issues of substantial computation load, large communication burden, user privacy disclosure, as well as limited scalability \cite{Yang2019review}.
Motivated by this, distributed optimization techniques, which allow computing nodes in a network to cooperatively optimize a global objective function with local communication only, are prosperously developed in recent years.
Distributed optimization over a networked system of $N$ nodes often takes the form of
\begin{equation}\label{eqn:min_sum_fi}
    \min_{x \in \Rbb^d} \sumi f_i(x),
\end{equation}
where the optimization variable $x \in \Rbb^d$ represents the global decision throughout the network, and $f_i(x)$ is the private convex objective function known only to node $i$. 
Applications of such a problem include but are not limited to power economic dispatch \cite{Yang2017EcoDisp,Yun2019EcoDisp}, distributed machine learning \cite{Nedic2020} and model predictive control \cite{Nedic2018}.

To date, there have been a variety of algorithms proposed to solve distributed optimization in the form of (\ref{eqn:min_sum_fi}) \cite{Nedic2009DGD,Shi2015PGEXTRA,Hong2017PGC,Aybat2018DPGA,Xu2018DFBBS,Makhdoumi2017DADMM}.
However, one common drawback of these methods is the lack of efficiency, i.e., a large number of iterations may be required to achieve a solution with desired accuracy, as the best convergence rate of many existing algorithms to solve the general convex form of (\ref{eqn:min_sum_fi}) is sublinear.
In industry, distributed optimization methods have not been widely adopted yet, since the substantial iterations ``increase the computation burden beyond the limit of practical interests" \cite{Wang2017} and result in degraded performance when the number of communication rounds is restricted.
Moreover, the worst-case convergence guarantees are provided by asymptotic analysis, while the practical performance of these optimization-based algorithms is sensitive to problem type, data pattern and parameter selection \cite{Balcan2021,De2021,Hutter2014}.
Thus, to achieve satisfactory problem-solving performance, a conventional approach is to manually try various algorithms to solve the problem of interest, hand-optimize the parameters, and then select the one with favorable performance, which is time-consuming and repetitive.

Recent advances in Learning to Optimize (L2O) have demonstrated the success of utilizing machine learning techniques to improve the performance of an existing optimization method or even learn a new optimization algorithm \cite{Chen2022L2Oreview}.
Examples include \textit{model-based} methods \cite{Chang2017,Meinhardt2017,Zhang2017,ZHANG2022180,Perdios2017} which fuse machine learning methods with analytical optimization algorithms, and \textit{model-free} methods \cite{Li2017,Bello2017,Andrychowicz2016,Chen2017,Cao2019} which develop completely new update rules by virtue of universal approximators such as neural networks.
However, most of the works are developed in \textit{centralized} configuration and do not specify the settings of distributed optimization.
In distributed optimization, the iterations are performed locally on each node through interactions among neighbors only and essentially need no global information. 
This update pattern is distinct from centralized methods and thus calls for an extended learning architecture.
Moreover, distributed optimization is essentially required to handle certain globally coupled constraints, which is not incorporated in typical centralized L2O frameworks.
For example, problem (\ref{eqn:min_sum_fi}) implicitly imposes a consensus constrain that requires the decisions of all the computing nodes to reach a consensus. 

Several L2O approaches under \textit{distributed} settings have been developed in recent years.
For example, learning-based ADMM variants are proposed in \cite{Graf2019,Biagioni2022,Zeng2022} to address practical distributed optimization problems such as wind farm control, coordination of demand response devices and optimal power flow.
Besides, \cite{Wang2021} unrolls two distributed algorithms Prox-DGD \cite{Zeng2018ProxDGD} and PG-EXTRA \cite{Shi2015PGEXTRA} with graph neural networks (GNNs) for decentralized statistical inference.
Despite the promising results, these L2O approaches mainly focus on improving the performance of a \textit{particularly selected} algorithm to solve a \textit{specific} problem, rather than a universal paradigm to solve the general form of problem (\ref{eqn:min_sum_fi}).
Therefore, the resulting problem-solving enhancement is restricted due to the lack of pre-selection among legitimate algorithms.
Moreover, the framework proposed in \cite{Wang2021} confines the number of running steps to at most the number of layers in the unrolled GNN, and the rest (e.g., \cite{Graf2019,Biagioni2022,Zeng2022}) require a centralized coordination scheme and are thus not fully decentralized.

In this paper, we establish a novel L2O framework to efficiently solve the \textit{general} form of problem (\ref{eqn:min_sum_fi}) and overcome the above drawbacks of the prior works.
We note that reinforcement learning (RL) has shown the ability to learn optimal policies in sequential decision-making problems, which is naturally suitable to model the sequential iterative process of numerical optimization.
Moreover, deep reinforcement learning (DRL) techniques utilize deep neural networks to extract features from high-dimensional inputs, so that problems with real-world complexity can be handled via end-to-end learning.
Inspired by these, we develop a DRL method to learn a policy for adaptive configuration within a unifying algorithmic form of distributed optimization.
This yields a model-based L2O framework that is (i) not confined to a single optimization algorithm, (ii) implemented in a fully decentralized manner, and (iii) highly efficient and flexible when applied in practice.

More specifically, we parameterize a surrogate-function-based unifying optimization method called Approximate Method of Multipliers \cite{Wu2022DAMM} as the \textit{base model}, which includes various first-order and second-order optimization algorithms along with the corresponding parameter ranges. 
Then, an RL agent is trained to make adaptive decisions regarding the algorithmic form of the base model for each iteration step according to the real-time information provided by computing nodes.
This model-based paradigm not only takes advantage of data-driven learning capacity, but also benefits from the convergence guarantee of the base model to solve (\ref{eqn:min_sum_fi}).
We exploit the aggregation of local consensus information, which measures the disagreements among neighboring nodes, and local objective information such as the gradients or Hessians of the local objective functions, to summarize the solving progress of a problem and to provide insights into the problem type as well as data pattern for the agent to make adaptive decisions.
This bundle of information constitutes the states observed by the RL agent.
We train this framework on a number of practical problem instances of similar structures yet different problem data using Proximal Policy Optimization (PPO) \cite{schulman2017PPO}, one of the most prevalent DRL methods owing to its sample efficiency and stable learning process.
In the experiments, we evaluate the learning-based optimization framework on various smooth and non-smooth classes of objective functions, and demonstrate that it outperforms several state-of-the-art distributed optimization algorithms in terms of convergence speed and solution accuracy.
Though the proposed framework is only trained over a relatively short time horizon, its superior performance is shown to be extended to subsequent iterations.

The rest of the paper is organized as follows:
Section \ref{sec:DRL} provides a preliminary overview on reinforcement learning, policy gradient methods, and PPO.
Section \ref{sec:prob_formu} formulates the general distributed convex optimization problem.
Section \ref{sec:method} elaborates the principles to build the learning-based framework for efficient distributed optimization.
Section \ref{sec:exps} demonstrates the effectiveness of the proposed method with various numerical examples.
Finally, Section \ref{sec:conclusion} concludes the paper.

\textbf{Notation:}
Let $A=(a_{i j})_{n \times n}$ be an $n \times n$ real matrix of which the $(i, j)$-entry is denoted by $[A]_{i j}$.
$\operatorname{Null}(A)$ denotes the null space of $A$.
$\mathbf{0}$ and $\mathbf{O}$ represent a zero vector and a zero square matrix of proper dimensions respectively.
$I_n$ denotes the $n \times n$ identity matrix. For two matrices $A$ and $B$, $A \otimes B$ represents the Kronecker product of $A$ and $B$.
We use $A \succeq B$ to represent that $A-B$ is positive semidefinite, if $A$ and $B$ are two square matrices of the same size.
For any vector $x \in \Rbb^n$ and matrix $A \succeq \mathbf{O}$, $\lVert x \rVert_1$ denotes the $\ell_1$-norm of $x$, $\lVert x \rVert=\sqrt{x^T x}$ and $\lVert x \rVert_A=\sqrt{x^T A x}$.
For any countable set $S$, $\lvert S \rvert$ denotes its cardinality.
We use $\mathcal{I}_X$ to denote the indicator function with respect to the set $X$, that is, $\mathcal{I}_X(x)=0$ if $x \in X$ and $\mathcal{I}_X(x)=+\infty$ otherwise.

\section{Preliminaries on Deep Reinforcement Learning}
\label{sec:DRL}
This section provides a brief background on DRL techniques.
Section \ref{subsec:MDP} introduces the Markov decision process (MDP), which is the model to formulate RL problems.
Section \ref{subsec:PG} discusses the classical policy gradient methods, and Section \ref{subsec:PPO} describes the proximal policy optimization (PPO) algorithm \cite{schulman2017PPO}, a DRL method originated from standard policy gradient approaches that is employed in this paper.

\subsection{Markov Decision Process} \label{subsec:MDP}

RL is a branch of machine learning techniques to address sequential decision problems.
It involves an agent interacting with an environment by observing the current state, taking an action according to a policy, receiving a reward, and transitioning to a new state.
The objective of RL is to find the optimal policy that maximizes the expected cumulative reward obtained by the agent over time.

Typically, an RL problem can be formulated as an MDP.
In this paper, we consider a finite-horizon MDP with continuous state and action spaces, which is characterized by the tuple $\mathcal{M}=(\mathcal{S}, \mathcal{A}, p_0, \mathcal{P}, R, \gamma)$, where $\mathcal{S}$ and $\mathcal{A}$ denote the state space and action space respectively, $p_0: \mathcal{S} \rightarrow \mathbb{R}^{+}$ is the probability density of initial states, $\mathcal{P}: \mathcal{S} \times \mathcal{A} \rightarrow \Delta_{\mathcal{S}}$ is the transition dynamics specifying a probability distribution over the successor state given the current state and action, $R: \mathcal{S} \times$ $\mathcal{A} \rightarrow \Rbb$ is the reward function, and $\gamma \in(0,1)$ is the discount factor that determines the relative importance of immediate rewards compared to future rewards.
A policy $\pi: \mathcal{S} \rightarrow \Delta_{\mathcal{A}}$ provides a mapping from a state to a probability distribution over actions at this state.
Particularly, we use $\pi (a \mid s)$ to denote the conditional probability of taking action $a$ given state $s$ according to the policy.

\subsection{Policy Gradient} \label{subsec:PG}

To secure the optimal policy in the an RL problem, policy gradient methods use samples collected from the agent-environment interaction to estimate the gradient of a performance measure $J(\theta)$ with respect to the policy parameters $\theta$, and then perform stochastic gradient ascent to update the parameterized policy $\pi_\theta$ in the direction of higher performance.
To reduce the variance of the gradient estimates and stabilize the learning process, policy gradient methods usually employ the actor-critic architecture, where the actor gradually adjusts the policy by updating the parameter $\theta$, and the critic estimates the value function $V(s)$ or action-value function $Q(s,a)$, which measures the expected reward accumulated from state $s$, for the actor to update the policy.
In real-world scenarios, deep neural networks are often utilized to approximate the policies and value functions, yielding various DRL techniques, which are able to handle high-dimensional inputs and learn complex, nonlinear relationships.

Policy gradient methods usually adopt the gradient estimation in the following form:
\begin{equation}\label{eqn:hat_g}
    \hat{g}=\hat{\mathbb{E}}_t\left[\nabla_\theta \log \pi_\theta\left(a_t \mid s_t\right) \hat{A}_t\right],
\end{equation}
where $\hat{A}_t$ is the estimation of the advantage function $A\left(s_t, a_t\right)=Q\left(s_t, a_t\right)-V\left(s_t\right)$ at timestep $t$, and $\hat{\mathbb{E}}_t[\cdot]$ denotes the empirical expectation over a finite batch of samples collected from the agent-environment interaction.
In practical implementation, the gradient estimation is often obtained by constructing a surrogate objective whose gradient is the same as (\ref{eqn:hat_g}), such as
\begin{equation}\label{eqn:L_PG}
    L^{PG}(\theta)=\hat{\mathbb{E}}_t\left[\log \pi_\theta\left(a_t \mid s_t\right) \hat{A}_t\right],
\end{equation}
and then using automatic differentiation softwares to differentiate this objective.

\subsection{Proximal Policy Optimization} \label{subsec:PPO}

The motivation behind PPO is to prevent the policy from drastic changes during the update process, since this may lead to unstable training or destructive updates.
PPO seeks to update the policy more conservatively by incorporating a clipping operator in the objective function to constrain each update, ensuring that the new policy remains close to the previous one.
Specifically, PPO aims to maximize the following surrogate objective function
$$
L^{CLIP}\!(\theta)\!=\!\hat{\mathbb{E}}_t\!\left[\min \left(r_t(\theta) \hat{A}_t, \operatorname{clip}\left(r_t(\theta), 1-\epsilon, 1+\epsilon\right) \hat{A}_t\right)\right]\!,
$$
where $r_t(\theta)$ denotes the probability ratio $r_t(\theta)=\frac{\pi_\theta\left(a_t, s_t\right)}{\pi_{\theta_{\text{old}}}\left(a_t, s_t\right)}$, the the $\operatorname{clip}(\cdot)$ operator saturates the probability ratio between $[1-\epsilon, 1+\epsilon]$.
For implementation, one simply constructs the objective $L^{CLIP}$ instead of $L^{PG}$ in (\ref{eqn:L_PG}), and uses an automatic differentiation software to perform stochastic gradient ascent on this objective for the policy updates.

\section{Problem Formulation}
\label{sec:prob_formu}
This section formulates the general distributed optimization problem that we investigate.

We consider a networked system modeled as an undirected, connected graph $\mathcal{G}=(\mathcal{V}, \mathcal{E})$, where $\mathcal{V}=\{ 1,2,\dots,N \}$ denotes the set of computing nodes and $\mathcal{E} \subseteq \left\{ \{i,j\} : i,j \in \mathcal{V}, i \neq j \right\}$ denotes the set of underlying communication links between neighboring nodes.
Each computing node $i\in\mathcal{V}$ possesses a private objective function $f_i : \Rbb^d \rightarrow \Rbb$, which is convex but not necessarily smooth, and exchanges selected information with its neighbors belonging to the set $\mathcal{N}_i := \left\{ j \in \mathcal{V}: \{i,j \in \mathcal{E}\} \right\}$.

Typically, all the computing nodes are required to address the following general distributed optimization problem in a composite form:
\begin{equation}\label{eqn:min_sum_si_p_ri}
    \min_{x\in\Rbb^d} \sumi \left( s_i(x)+r_i(x) \right),
\end{equation}
where each $s_i : \Rbb^d \rightarrow \Rbb$ is a convex function and has a Lipschitz continuous gradient, and each $r_i : \Rbb^d \rightarrow \Rbb \cup \{+\infty\}$ is a convex and possibly non-differentiable function.
Both $s_i$ and $r_i$ are allowed to be a zero function, yielding smooth or non-smooth objective functions.
Notice that $r_i$ can also contain an indicator function $\mathcal{I}_{X_i}$ with respect to a closed convex set $X_i \subset \Rbb^d$, so that non-smooth, constrained convex programs can also be formulated by (\ref{eqn:min_sum_si_p_ri}).

Many practical engineering problems can be cast into the form of (\ref{eqn:min_sum_si_p_ri}).
In distributed learning \cite{Predd2006}, $s_i$ can be the mean squared error (MSE) or cross-entropy loss, while $r_i$ can represent an $\ell_1$-norm regularizer or the mean absolute error (MAE) loss.
Also, many real-world constrained optimization problems such as coordination of distributed energy resources (DERs) \cite{Yang2019review} can be formulated as (\ref{eqn:min_sum_si_p_ri}) via dualization.

Let each computing node $i\in\mathcal{V}$ maintain a local decision $x_i \in \Rbb^d$.
Using the global compact variable $\xb=(x_1^T,\dots,x_N^T)^T \in \Rbb^{Nd}$, we can rewrite (\ref{eqn:min_sum_si_p_ri}) as the following constrained problem:
\begin{equation} \label{eqn:min_s_p_r}
\begin{array}{cl}
    \min\limits_{\xb\in\Rbb^{Nd}} & s(\xb)+r(\xb) \\
    \st & H^{\frac{1}{2}}\xb = \mathbf{0},
\end{array}
\end{equation}
where $s(\xb):=\sumi s_i(x_i)$, $r(\xb):=\sumi r_i(x_i)$ and $H \in \Rbb^{Nd\times Nd}$ is such that $H=H^T$, $H\succeq\mathbf{O}$, $\operatorname{Null}(H)=S$.
It is shown in \cite{Mokhtari2016} that the constraint $H^{\frac{1}{2}}\xb = \mathbf{0}$ in (\ref{eqn:min_s_p_r}) is identical to the consensus constraint $\xb \in S:=\left\{\xb \in \Rbb^{N d} \mid x_1=\cdots=x_N\right\}$, so that $H^{\frac{1}{2}}\xb$ quantifies the disagreements among all the computing nodes.

\section{Deep Reinforcement Learning Approach}
\label{sec:method}
To solve problem (\ref{eqn:min_s_p_r}) efficiently and reliably, we establish a model-based L2O framework, which benefits from both the convergence guarantee of the base model and the high efficiency resulted from data-driven learning.
We introduce the base model to solve the problem (\ref{eqn:min_s_p_r}) in Section \ref{subsec:base_model}, and elaborate how this base model is integrated with DRL techniques in Section \ref{subsec:int_with_DRL}.

\subsection{Base Model} \label{subsec:base_model}

We consider an optimization-theoretic base model, which ensures the convergence performance while solving problem (\ref{eqn:min_s_p_r}).
Notice that \cite{Wu2022DAMM} has developed an optimization method, called Approximate Method of Multipliers (AMM), supporting the general formulation (\ref{eqn:min_s_p_r}) and capable of achieving the convergence rate of best order $O(1/k)$.
AMM utilizes a surrogate function to approximate the behavior of the Method of Multipliers \cite{Wu2022DAMM}.
By choosing different forms of the surrogate function, AMM unifies a variety of existing and new distributed optimization algorithms.
This motivates us to investigate the autonomous, adaptive configuration within the unifying algorithmic form of AMM and later integrate it with DRL techniques to achieve more efficient problem solving.

Specifically, to solve (\ref{eqn:min_s_p_r}), AMM starts from initializing the dual variable $\vb^0 \in \Rbb^{Nd}$, and updates according to
\begin{align}
    & \xb^{k+1}\!=\!\argminxb u^k(\xb)+r(\xb)+\frac{\rho}{2}\lVert\xb\rVert_H^2+(\vb^k)^T H^{\frac{1}{2}} \xb, \label{eqn:AMM_x}\\
    & \vb^{k+1}\!=\!\vb^k+\rho H^{\frac{1}{2}} \xb^{k+1}, \quad \forall k \geq 0, \label{eqn:AMM_v}
\end{align}
where $\rho>0$ is the penalty parameter, and $u^k:\Rbb^{Nd}\rightarrow\Rbb$ is a convex and smooth surrogate function at iteration $k$ to approximate the Method of Multipliers such that $\nabla u^k(\xb^k)=\nabla s(\xb^k)$ and the function $u^k + \frac{\rho}{2} \lVert\cdot\rVert_H^2$ is strongly convex \cite[Assumption 2]{Wu2022DAMM}.

To perform the updates in a decentralized manner, we adopt a distributed implementation of AMM in \cite{Wu2022DAMM} capable of handling the composite form (\ref{eqn:min_s_p_r}), which operates as follows:
Each node $i\in\mathcal{V}$ arbitrarily initializes its local primal variable $x_i^0 \in \mathbb{R}^d$ and dual variable $q_i^0 \in \mathbb{R}^d$ such that $\sum_{i\in\mathcal{V}} q_i^0 = \mathbf{0}$, and then all nodes in the network $\mathcal{G}$ iteratively update as
\begin{align}
    &x_i^{k+1} = \argminx \psi_i^k(x) + r_i(x) + \langle \tilde{c}_i^k, x \rangle, &\forall k \geq 0, \label{eqn:DAMM_x}\\
    &q_i^{k+1} = q_i^k + \rho\sum_{j\in\mathcal{N}_i\cup\{i\}}p_{ij}x_j^{k+1}, &\forall k \geq 0, \label{eqn:DAMM_q}
\end{align}
where $\tilde{c}_i^k := q_i^k-\nabla\psi_i^k(x_i^k)+\nabla s_i(x_i^k)+\rho\sum_{j\in\mathcal{N}_i\cup\{i\}}p_{ij}x_j^k$, $\forall k \geq 0$, $p_{ij} = p_{ji}<0$, $\forall\{i,j\}\in\mathcal{E}$, $p_{ii}=-\sum_{j\in\mathcal{N}_i} p_{ij}$, $\forall i \in \mathcal{V}$, and $\psi_i^k(\cdot)$ are smooth and strongly convex functions such that $\left(\sumi \psi_i^k\left(x_i\right)\right)-\frac{\rho}{2}\lVert \xb \rVert_H^2$ are convex.
All of the weights $p_{ij}$ constitute a weight matrix $P$ such that $[P]_{ij}=p_{ij}\ \forall i,j \in \mathcal{V}$ and $H=P\otimes I_d$.
Indeed, $\psi_i^k(\cdot)$ are local, possibly time-varying functions that can be personalized by each node.

To build the base model, we allow the function $\psi_i^k(\cdot)$ to include the \textit{second-order} information of the local smooth component $s_i(\cdot)$, so that the resulting algorithm can be first-order, second-order or hybrid.
We then parameterize $\psi_i^k(x)$ as the following quadratic form:
\begin{equation}\label{eqn:param_psi}
    \psi_i^k(x)=\frac{1}{2}x^T\left(\alphaik\nabla^2s_i(x_i^k)+\betaik I_d\right)x,
\end{equation}
where $\alphaik, \betaik \geq 0$ are possibly time-varying parameters that can be determined by each node individually.
Additionally, we introduce more flexibility to the penalty parameter $\rho$ by casting it to a similar personalized and time-varying form, i.e., each node $i$ maintains a private penalty parameter $\rho_i^k$.
Then, our base model is given by
\begin{align}
    &x_i^{k+1}\! =\! \argminx\! \frac{1}{2}\lVert x \rVert^2_{\alphaik\nabla^2s_i(x_i^k)+\betaik I_d}\! +\! r_i(x)\! + \!\langle c_i^k, x \rangle, \label{eqn:param_DAMM_x} \\
    &q_i^{k+1} = q_i^k + \rho_i^k\sum_{j\in\mathcal{N}_i\cup\{i\}}p_{ij}x_j^{k+1}, \label{eqn:param_DAMM_q}
\end{align}
where $c_i^k = q_i^k-\left(\alphaik\nabla^2s_i(x_i^k)+\betaik I_d\right)x_i^k+\nabla s_i(x_i^k)+\rho_i^k\sum_{j\in\mathcal{N}_i\cup\{i\}}p_{ij}x_j^k$.

Note that with all $\alpha_i^k>0$, (\ref{eqn:param_DAMM_x})\textminus(\ref{eqn:param_DAMM_q}) form a distributed \textit{second-order} method for solving problem (\ref{eqn:min_s_p_r}), which inherits the sublinear convergence rate of AMM if $\sumi \left( \psi_i^k(x_i)-\frac{\rho_i^k}{2}\lVert x_i \rVert_P^2 \right)$ is convex.
If we let all $\alpha_i^k=0$, the function $\psi_i^k$ reduces to $\psi_i^k(\cdot) = \frac{\beta_i^k}{2} \lVert \cdot \rVert^2$, so that (\ref{eqn:param_DAMM_x})\textminus(\ref{eqn:param_DAMM_q}) become distributed first-order methods, which still generalize several existing \textit{first-order} algorithms such as PGC \cite{Hong2017PGC}, PG-EXTRA \cite{Shi2015PGEXTRA}, DPGA \cite{Aybat2018DPGA}, and D-FBBS \cite{Xu2018DFBBS}, as is shown in \cite{Wu2022DAMM}.
Therefore, the base model given by (\ref{eqn:param_DAMM_x})\textminus(\ref{eqn:param_DAMM_q}) is naturally endowed with the generality to unify various algorithm types.

\subsection{Integration with DRL} \label{subsec:int_with_DRL}

While the base model (\ref{eqn:param_DAMM_x})\textminus(\ref{eqn:param_DAMM_q}) inherits the worst-case convergence guarantee of AMM, its practical performance still suffers from the sensitivity to problem types, data pattern and parameter choices.
To overcome this issue, below we develop a DRL method to learn a policy for autonomous and adaptive configuration of the base model parameterized via (\ref{eqn:param_psi}) to perform each local iteration.

\begin{figure}[t]
    \centering
    \includegraphics[width=0.46\linewidth]{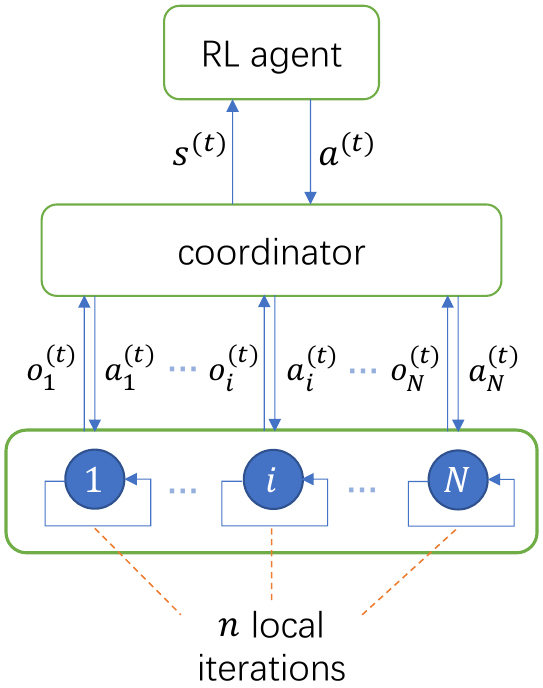}
    \caption{Interaction paradigm of the learning-based framework within a communication round. Circles marked with $i=1,\dots,N$ represent the computing nodes in the networked system.}
    \label{fig:framework}
\end{figure}

To properly formulate the RL problem, we propose the learning-based framework as shown in Figure \ref{fig:framework}, which includes an RL agent interacting with the computing nodes via a \textit{coordinator} that is designed with \textit{no computing ability}.
In this framework, the RL environment is regarded as the solving process of problem (\ref{eqn:min_s_p_r}) under the base model (\ref{eqn:param_DAMM_x})\textminus(\ref{eqn:param_DAMM_q}).
The interaction between the RL agent and the environment is described as follows:
In each communication round $t=0,1,\dots$, the RL agent observes the current solving progress of the problem, represented by the state $s^{(t)}$, and determines the local update rules on computing nodes in this round according to the current policy.
Such a decision is represented by the action $a^{(t)}$.
Then, all the computing nodes in the network operate as the RL agent assigns, and feed a reward $R(s^{(t)}, a^{(t)})$ back to the agent, measuring how beneficial this action is to accelerate the solving progress.
The above interactive process within a communication round constitutes a timestep in the typical RL context for training the agent.

In our proposed framework, the communication between the computing nodes and the RL agent is bridged via the coordinator.
In each communication round, the coordinator aggregates the local observation $o_i^{(t)}$ of each node $i\in\mathcal{V}$ to constitute the state $s^{(t)}$.
This state information is then directed to the RL agent.
Once the RL agent determines the global action $a^{(t)}$ for the whole network according to $s^{(t)}$, the coordinator receives $a^{(t)}$ from the RL agent and then dispatches the corresponding local action $a_i^{(t)}$ to each computing node.
Subsequently, each node in the network starts its local update, and the coordinator will \textit{not} participate in this distributed computing process.
In our design, the local update of each node in a round is allowed to involve \textit{multiple} steps of the iteration (\ref{eqn:param_DAMM_x})\textminus(\ref{eqn:param_DAMM_q}), obviating the need for frequent communication between the computing nodes and the coordinator.
Hence, the relation between the communication round index $t$ and the iteration step index $k$ in Section \ref{subsec:base_model} is $k=nt$, where $n>0$ is the number of local iterations (\ref{eqn:param_DAMM_x})\textminus(\ref{eqn:param_DAMM_q}) in each round.

\begin{remark}
The proposed framework is essentially different from the existing methods of centralized optimization or distributed optimization with a central node.
In these methods, the central processing platform either performs the global updates directly or is indispensable for every local update of the computing nodes. 
Thus, its failure will result in the collapse of the entire networked system.
Our method, on the contrary, is more robust since the coordinator does not participate in the optimization process.
Even if the coordinator fails at some time, the whole system will continue operating according to the previous action, and eventually solve the optimization problem (in a less efficient way) due to the convergence guarantee of the base model.
\end{remark}

In the following parts, we describe the settings of the state space, action space and reward function in detail, which are crucial to achieve an effective and stable training process.

\textbf{State space:}
The state should contain necessary information to characterize the solving progress of the distributed optimization problem (\ref{eqn:min_s_p_r}).
In centralized optimization, it is natural to use first-order information such as the gradients at the iteration points to trace the solving progress, as they gradually approach zero while the optimization algorithms proceed.
However, this is not applicable to distributed optimization since the local gradients provide no conclusive information on the global optimality.

To reveal the solving progress of problem (\ref{eqn:min_s_p_r}), first note that the optimality conditions for (\ref{eqn:min_s_p_r}) are equivalent to the primal and dual feasibility as follows:
$$
    H^{\frac{1}{2}} \xb^\star=\mathbf{0}, \quad
    -\nabla s(\xb^\star) - H^{\frac{1}{2}} \vb^\star \in \partial r(\xb^\star),
$$
where $(\xb^\star, \vb^\star)$ is the primal-dual optimum pair.
In addition, (\ref{eqn:AMM_x}) in the AMM update rule is equivalent to obtaining the unique $\xb^{k+1}$ such that
\begin{equation}
    -\nabla u^k(\xb^{k+1}) - \rho H \xb^{k+1} - H^{\frac{1}{2}} \vb^k \in \partial r(\xb^{k+1}).
\end{equation}
This, together with (\ref{eqn:AMM_v}) and $\nabla u^k(\xb^k)=\nabla s(\xb^k) \ \forall k \geq 0$ in the AMM assumptions, leads to 
\begin{equation}
    -\nabla s(\xb^{k+1}) - H^{\frac{1}{2}} \vb^{k+1} \in \partial r(\xb^{k+1}),
\end{equation}
which means that $(\xb^{k+1}, \vb^{k+1})$ is always dual feasible during the process of AMM.
Then the primal feasibility $H^{\frac{1}{2}}\xb^k = \mathbf{0}$ is sufficient to indicate the optimality when we use AMM to solve problem (\ref{eqn:min_s_p_r}).
On the other hand, as is shown in \cite{Wu2022DAMM}, the primal residual $H^{\frac{1}{2}}\xb^k$, which measures the disagreements among all the computing nodes, converges to zero as AMM proceeds.
Hence, the consensus information contained in $H^{\frac{1}{2}}\xb^k$ can be exploited by the RL agent to represent the solving progress of (\ref{eqn:min_s_p_r}).

However, the global consensus information $H^{\frac{1}{2}}\xb^k$ is inaccessible to the computing nodes in the designed distributed framework.
Additionally, the coordinator cannot directly collect local decision variables $x_1^k, \dots, x_N^k$ to compute the global consensus $H^{\frac{1}{2}}\xb^k=(P\otimes I_d)^{\frac{1}{2}}\xb^k$.
There are several reasons: 
(1) additional global coordination is required to obtain the weight matrix $P$, (2) collecting all local decision variables increases the risk of privacy disclosure, and (3) the coordinator is designed to have no computing ability, so that all the computations can only be performed locally on the computing nodes.
Fortunately, the base model (\ref{eqn:param_DAMM_x})\textminus(\ref{eqn:param_DAMM_q}) enables each node to individually obtain the local consensus information $\sum_{j\in\mathcal{N}_i\cup\{i\}}p_{ij}x_j^k$ via communicating with its neighbors.
Therefore, we can use the aggregation of this local consensus information from all nodes, instead of $H^{\frac{1}{2}}\xb^k$, to reflect the state of convergence.

Moreover, note that if the objective function has a smooth component (i.e., $s_i(\cdot)\neq 0$), the base model (\ref{eqn:param_DAMM_x})\textminus(\ref{eqn:param_DAMM_q}) uses first-order and second-order information of the smooth component to update.
Hence, we naturally expect such information to guide the adaptive configuration of the base model, and include this information in the state.
Here, the second-order information is extracted through the eigenvalues of $\nabla^2 s_i(x_i^k)$, which not only reduces the dimensions required to store the second-order information from $d \times d$ to $d$, but also reveals necessary geometric regularities of the objective function as well as the data pattern \cite{sagun2017eigenvalues}.

The above local information with $n$-iteration (the past round) history constitutes the local observation $o_i^{(t)}$ of each node $i\in\mathcal{V}$ for each communication round $t$:
$$o_i^{(t)}\! = \!\!\begin{bmatrix} 
    (\sigma_i^{n(t-1)+1}, g_i^{n(t-1)+1}, \lambda_i^{n(t-1)+1}),\! \cdots\!, (\sigma_i^{nt}, g_i^{nt}, \lambda_i^{nt})
\end{bmatrix}\!,$$
where $\sigma_i^k = \sum_{j\in\mathcal{N}_i\cup\{i\}}p_{ij}x_j^{k}$, $g_i^k = \nabla s_i(x_i^{k})$, and $\lambda_i^k \in \Rbb^d$ consists of the eigenvalues of the Hessian $\nabla^2 s_i(x_i^k)$, $\forall k\geq0$.
The state $s^{(t)}$ is the aggregation of the local observations:
$$s^{(t)} = \begin{bmatrix}
    o_1^{(t)},\cdots,o_N^{(t)}
\end{bmatrix}.$$
Intuitively, the state $s^{(t)}$ provides regularities of the problem instance such as the problem type and data distribution, and summarizes the solving progress of this problem.

\textbf{Action space:}
In communication round $t$, the (global) action $a^{(t)}$ specifies the exact form of the parameterized base model (\ref{eqn:param_DAMM_x})\textminus(\ref{eqn:param_DAMM_q}) to perform local updates.
To control all the local iterations in a round, a tuple of three parameters $(\alphaik, \betaik, \rhoik)$ at iteration steps $k=n(t-1)+1,\cdots,nt$ should be dispatched to each node $i\in\mathcal{V}$.
This gives rise to an action space with large dimensions.
For simplicity, we adopt the following setting:
In round $t$ the RL agent determines only one tuple of the parameters $(\alpha^{(t)}, \beta^{(t)}, \rho^{(t)})$, and the coordinator then dispatches these parameters to all the computing nodes for every local iteration in the round.
Therefore, the action of the RL agent has the form
$$a^{(t)} = (\alpha^{(t)}, \beta^{(t)}, \rho^{(t)})$$
with $\alpha^{(t)}, \beta^{(t)} \geq 0$ and $\rho^{(t)} > 0$, while the coordinator dispatches the local actions as $a_1^{(t)}=\cdots=a_N^{(t)}=a^{(t)}$.
In practice, $\alpha^{(t)}, \beta^{(t)}, \rho^{(t)}$ can be bounded from above by  some customized constants to stabilize the learning process.
We note, however, that it is straightforward to extend this global action to a node-specific paradigm when personalized local actions for each computing node are necessary:
$$a^{(t)} = \begin{bmatrix}
    (\alpha_1^{(t)}, \beta_1^{(t)}, \rho_1^{(t)}),\cdots,(\alpha_N^{(t)}, \beta_N^{(t)}, \rho_N^{(t)})
\end{bmatrix}.$$

\textbf{Reward function:}
As our goal is to accelerate the solving processes of distributed optimization problems, the reward function should be designed to effectively penalize policies that lead to slow convergence.
In distributed optimization (\ref{eqn:min_s_p_r}), the performance metrics typically include the objective error $\lvert \sumi (s_i(x_i^k)+r_i(x_i^k)-s_i(x^\star)-r(x^\star)) \rvert$ and consensus error $\sumi \lVert x_i^k - \bar{x}^k \rVert^2$, where $x_i^k$ is the local iterate of node $i$ at iteration $k$, $x^\star$ is the ground-truth decision variable, and $\bar{x}^k=\frac{1}{N}\sumi x_i^k$ is the mean of all local iterates at iteration $k$.
These metrics can be unified by the mean square error (MSE) of all local iterates:
$$\operatorname{MSE}(\xb^k, x^\star) = \frac{1}{N} \sumi \lVert x_i^k - x^\star \rVert^2.$$
Therefore, we define the reward function as the total MSE accumulated from the round $t$:
$$R(s^{(t)}, a^{(t)}) = \sum_{k=nt+1}^{n(t+1)} \operatorname{MSE}(\xb^k, x^\star).$$
This reward setting encourages all the computing nodes to reach the optimal decision variable $x^\star$ as quickly as possible.
In the training stage, the ground-truth $x^\star$ can be acquired via any commercial solver for learning the policy.
After the training, the ground-truth $x^\star$ is no longer needed for the learned framework to optimize other objective functions with similar structures.

With the above specific settings of the state space, action space and reward function, we train the RL agent in Figure \ref{fig:framework} using PPO, as described in Section \ref{subsec:PPO}, so that it learns the optimal policy for autonomous and adaptive configuration of the base model to perform local iterations on the computing nodes.

\section{Numerical Experiments}
\label{sec:exps}
In this section, we demonstrate the competence of the proposed learning-based framework to address the general distributed convex optimization (\ref{eqn:min_sum_si_p_ri}) via numerical experiments.
The experiments involve three problem cases with smooth or non-smooth objectives that correspond to the loss functions of regression or classification tasks.
Specifically, we consider: (A) linear least square regression with Lasso regularization \cite{Tibshirani1996}, in which the objective is in the composite form as (\ref{eqn:min_sum_si_p_ri}); (B) logistic regression, in which the objective is a special case of (\ref{eqn:min_sum_si_p_ri}) with all $r_i(\cdot)=0$; and (C) linear $\ell_1$-regression with Lasso regularization, in which the objective is another special case of (\ref{eqn:min_sum_si_p_ri}) with all $s_i(\cdot)=0$.
These three cases together provide a comprehensive illustration of our proposed framework to handle practical distributed optimization problems.

In the following experiments, we consider a networked system consisting of $\lvert \mathcal{V} \rvert = N = 10$ computing nodes and $\lvert \mathcal{E} \rvert = 30$ undirected communication links, structured with a coordinator and an RL agent as in Figure \ref{fig:framework} to perform learning-based optimization.
Each computing node agrees with its neighbors on setting $p_{ij}=p_{ji}=-\frac{1}{\max \left\{\left|\mathcal{N}_i\right|,\left|\mathcal{N}_j\right|\right\}+1}$, $\forall\{i, j\} \in \mathcal{E}$.
The number of local iterations in each communication round is set to $n=10$.
Initially, each node in the networked system performs $n$ local iterations with arbitrary local actions to form the initial state $s^{(0)}$ observed by the RL agent.

We employ the actor-critic architecture to train the RL agent. 
Specifically, we use a shallow neural network with two fully-connected hidden layers to model the mean of the policy.
A Gaussian head is attached to this network to learn the state-independent variance \cite{pfrl}.
With every input of the state, the policy network outputs the mean and variance of a Gaussian distribution where the action is sampled from.
We use another similar neural network with two fully-connected hidden layers to model the value function.

For each of the problem cases, we construct a dataset that contains 100 instances for training, 10 instances for validation, and 10 instances for test.
The instances of a problem case correspond to objective functions with the similar structure but different problem data, which will be specified later.
Then the RL agent is trained with PPO for 10 timesteps, which correspond to 10 communication rounds of the networked system, i.e., a total of $10\times 10=100$ local iterations (\ref{eqn:param_DAMM_x})\textminus(\ref{eqn:param_DAMM_q}) for each computing node.

To accelerate the training progress, we initially pretrain the mean of the policy to match a \textit{baseline} action using supervised learning.
The baseline can be any hand-selected action that simply leads to \textit{not inferior} performance of the base model (\ref{eqn:param_DAMM_x})\textminus(\ref{eqn:param_DAMM_q}).
For example, we may arbitrarily select the fixed action, say, $a = (\alpha^{(t)}, \beta^{(t)}, \rho^{(t)}) = (5,5,5)$ for all the communication rounds $t=1,\dots,10$, and check the base model under this fixed action on any instance(s) of the training set. 
If we find out that this action yields a descending trend of $\operatorname{MSE}(\xb^k, x^\star)$, then it can be selected as the baseline.
This pretrained mean along with the initialized variance constitutes the \textit{initial policy} before training, eliminating the tedious random exploration in the former stage of the training process.

\subsection{Linear Least Square Regression with Lasso Regularization} \label{subsec:exp_composite}

In this subsection, we consider learning a linear model by minimizing the Lasso regularized MSE loss over the aforementioned networked system:
\begin{equation}\label{prob:least_square_lasso_reg}
    \min_{x\in\mathbb{R}^d} \sum_{i=1}^{N} \left( \frac{1}{m_i} \sum_{j=1}^{m_i} \frac{1}{2}\left(a_{i,j}^T x - b_{i,j}\right)^2 + \lambda \lVert x \rVert_1 \right),
\end{equation}
where $a_{i,j} \in \Rbb^d$ and $b_{i,j}\in\Rbb$ denote the feature vector and label of the $j$-th sample on computing node $i$ respectively, $m_i$ is the number of data samples $\{a_{i,j}, b_{i,j}\}$ on each node, and $\lambda > 0$ is the regularization parameter.
The objective function has the composite form as (\ref{eqn:min_sum_si_p_ri}) with $s_i(x)=\frac{1}{m_i} \sum_{j=1}^{m_i} \frac{1}{2}\left(a_{i,j}^T x - b_{i,j}\right)^2$ and $r_i(x)=\lambda \lVert x \rVert_1$.

The experiment is conducted on a standard real dataset \textit{abalone} from UCI Machine Learning Repository \cite{abalone:1995}.
Each problem instance is generated by sampling $\sumi m_i=100$ feature-label pairs $\{a_{i,j}, b_{i,j}\}$ randomly and allocating them to all the computing nodes.
Therefore, each instance in the training set, validation set or test set corresponds to the same problem (\ref{prob:least_square_lasso_reg}) yet with different problem data.
The ground-truth decision variable $x^\star$ for each problem instance is obtained via solving (\ref{prob:least_square_lasso_reg}) in a centralized manner using CVXPY \cite{diamond2016cvxpy}.
We use the training set to train the RL agent and evaluate the performance of the current policy on the validation set during the training process.
Then we choose the policy which performs the best on the validation set, and report the average performance of this policy on the test set that the RL agent has never seen before.

\begin{figure}[t]
    \centering
    \begin{subfigure}[b]{0.99\linewidth}
      \includegraphics[width=\linewidth]{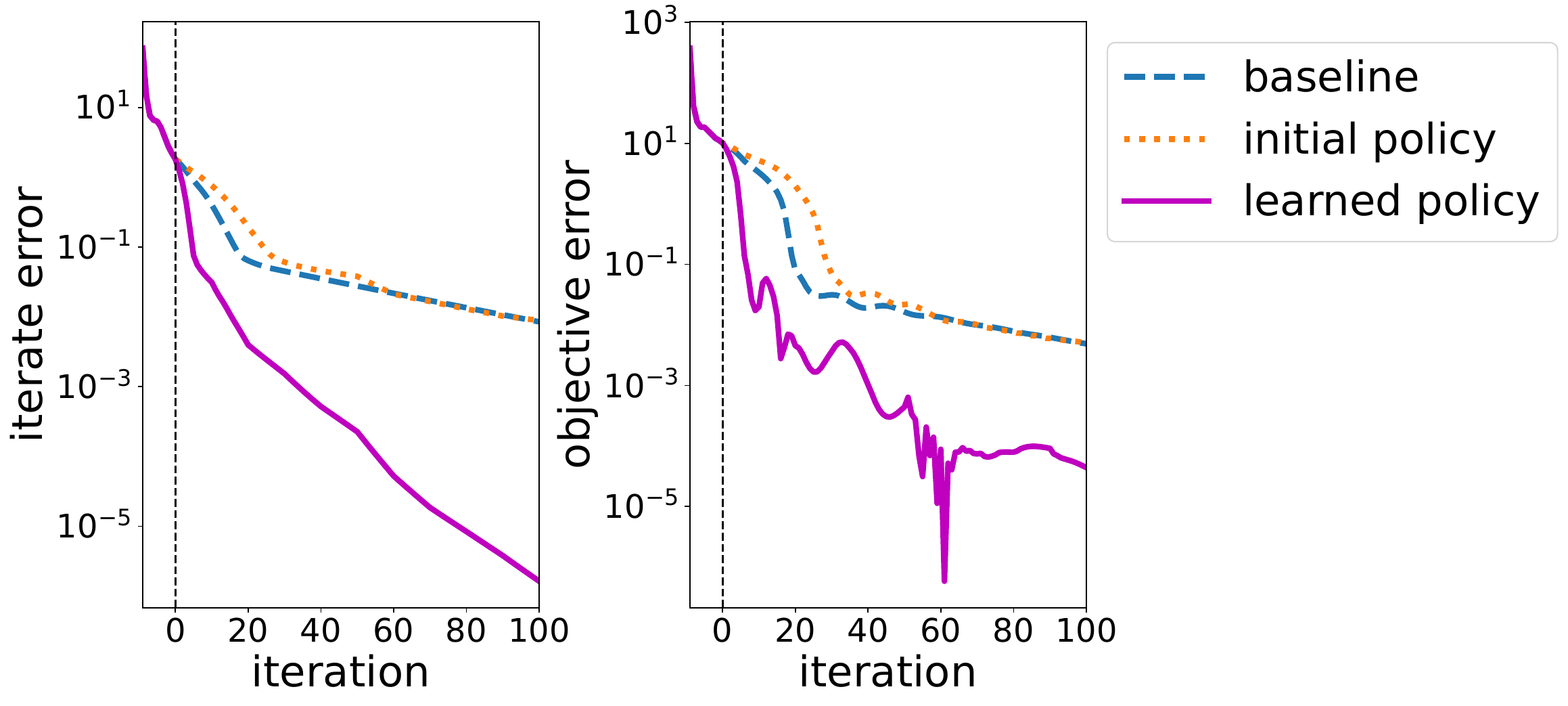}
      \caption{}
      \label{subfig:composite:train}
    \end{subfigure}
    \begin{subfigure}[b]{0.98\linewidth}
      \includegraphics[width=\linewidth]{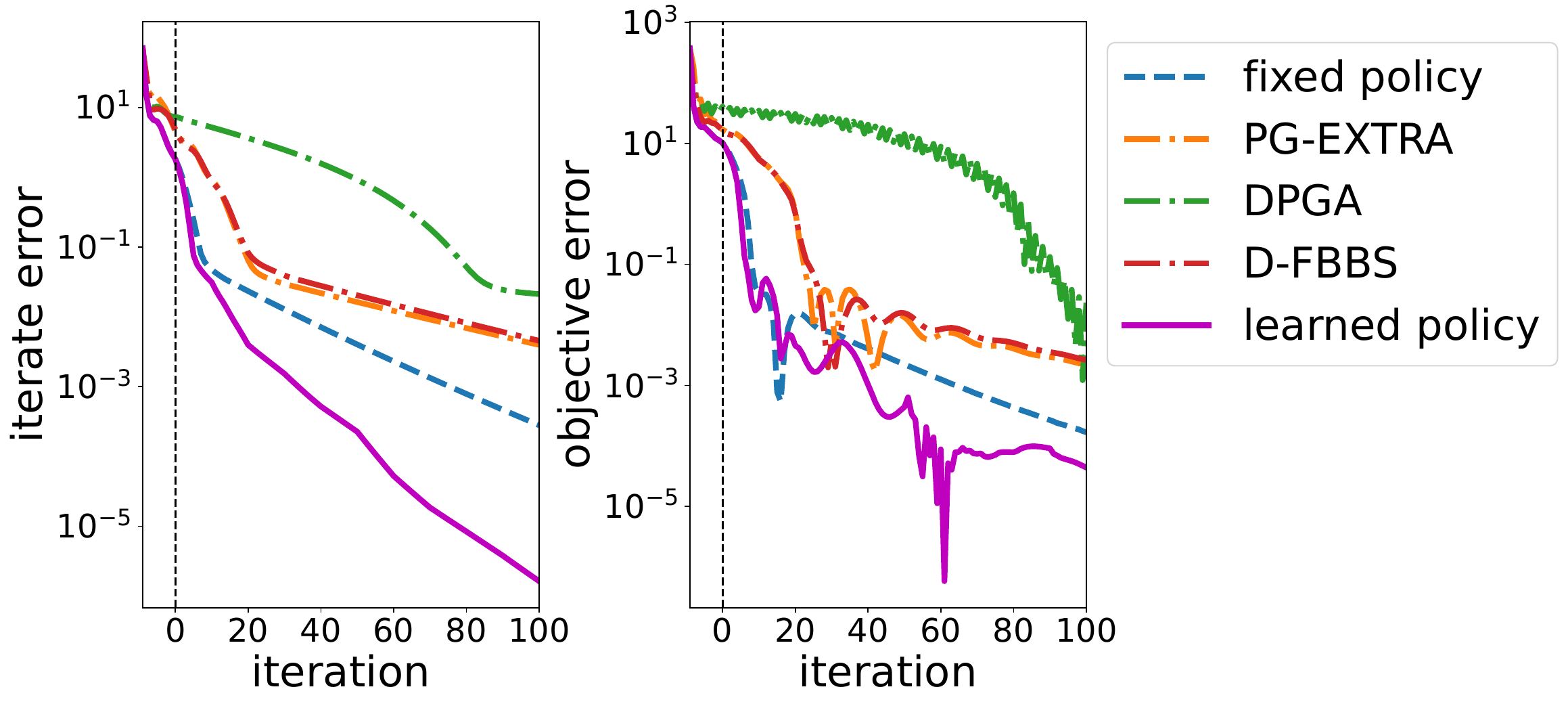}
      \caption{}
      \label{subfig:composite:compare}
    \end{subfigure}
    \caption{(a) Convergence performance of base model (\ref{eqn:param_DAMM_x})\textminus(\ref{eqn:param_DAMM_q}) under the baseline, the initial policy and the learned policy for solving (\ref{prob:least_square_lasso_reg}). (b) Convergence performance of base model (\ref{eqn:param_DAMM_x})\textminus(\ref{eqn:param_DAMM_q}) under the learned policy and the fixed policy (i.e., $\pi(a^c \mid s)\equiv 1$), as well as the convergence performance of state-of-the-art algorithms applicable to (\ref{prob:least_square_lasso_reg}).}
    \label{fig:composite}
\end{figure}

Figure \ref{fig:composite} demonstrates the effectiveness of our learning-based optimization framework by illustrating the learning process and comparing its convergence performance with state-of-the-art distributed optimization algorithms.
The iteration steps before the zero coordinate (vertical dashed line) are used to form the initial state $s^{(0)}$ for the learning-based framework.
In Figures \ref{subfig:composite:train}, we plot the convergence performance of the networked system under the baseline, initial policy and the learned policy in terms of the iterate error $\frac{1}{N} \sumi \lVert x_i^k - x^\star \rVert^2$ and the objective error $\lvert \sumi (s_i(x_i^k)+r_i(x_i^k)-s_i(x^\star)-r(x^\star)) \rvert$ at each iteration $k=0,\dots,100$.
The initial policy is intended to match the baseline with a variance for exploration.
As shown in Figure \ref{subfig:composite:train}, even though neither the baseline nor the initial policy yields satisfactory performance at first, the learning-based framework is able to achieve much faster convergence speed and solutions with much higher accuracy after training.

In Figures \ref{subfig:composite:compare}, we compare the learned policy with a \textit{fixed policy}, which controls the base model (\ref{eqn:param_DAMM_x})\textminus(\ref{eqn:param_DAMM_q}) using a hand-selected constant action $a^c$ (i.e., $\pi(a^c \mid s)\equiv 1$), via evaluating the convergence performance of the base model under these policies respectively.
To further indicate the superior performance of our learning-based framework, we incorporate various state-of-the-art distributed optimization methods into the comparison, which are capable of solving the general form of problem (\ref{eqn:min_sum_si_p_ri}) with the best $O(1/k)$ convergence rate, including PG-EXTRA \cite{Shi2015PGEXTRA}, DPGA \cite{Aybat2018DPGA} and D-FBBS \cite{Xu2018DFBBS}.
For fair comparison, these algorithms start with the same initialization as the learning-based method, and since they do not require the initial $n=10$ iteration steps to form the state, a total of 110 iteration steps are actually performed.
The constant action of the fixed policy and other methods for comparison are all fine-tuned using validation instances.
Apparently, the learned policy shows better convergence performance compared with the fixed policy, indicating that our learning-based framework surpasses the traditional hand-selecting paradigm for algorithm tuning.
In addition, base model (\ref{eqn:param_DAMM_x})\textminus(\ref{eqn:param_DAMM_q}) under the learned policy is shown to outperform PG-EXTRA, DPGA and D-FBBS in convergence speed and solution accuracy to a great extent.

\subsection{Logistic Regression} \label{subsec:exp_smooth}

Next, we consider learning a linear classifier by minimizing the $\ell_2$ regularized logistic loss over the networked system:
\begin{equation}\label{prob:logistic_reg}
  \min_{x\in\mathbb{R}^d} \sum_{i=1}^{N}\! \left(\! \frac{1}{m_i}\! \sum_{j=1}^{m_i}\! \left( -b_{i,j}a_{i,j}^T x \!+\! \log\left(\!1\!+\!e^{a_{i,j}^T x}\!\right)\! \right)\! +\! \frac{\lambda}{2} \lVert x \rVert_2^2\! \right)\!,
\end{equation}
where $a_{i,j} \in \Rbb^d$ and $b_{i,j}\in\{0,1\}$ denote the feature vector and category label of the $j$-th sample on computing node $i$ respectively, $m_i$ denotes the number of data samples $\{a_{i,j}, b_{i,j}\}$ on each node, and $\lambda > 0$ is the regularization parameter.
This smooth problem is the special case of (\ref{eqn:min_sum_si_p_ri}) with $s_i(x)=\frac{1}{m_i} \sum_{j=1}^{m_i} \left( -b_{i,j}a_{i,j}^T x + \log\left(1+e^{a_{i,j}^T x}\right) \right) + \frac{\lambda}{2} \lVert x \rVert_2^2$ but all $r_i(x)=0$.

The experiment is conducted on another real dataset \textit{breast cancer Wisconsin} from UCI Machine Learning Repository \cite{breast_cancer_wisconsin_original:1992} for this logistic regression problem.
The problem instances in the training set, validation set and test set are generated following the same procedures as described in Section \ref{subsec:exp_composite}, such that each instance corresponds to problem (\ref{prob:logistic_reg}) with different problem data.
Similarly, we train the RL agent on the training set, select the best policy according to the validation loss, and report the average performance of this policy on the test set that is new to the RL agent.

\begin{figure}[t]
    \centering
    \begin{subfigure}[b]{0.99\linewidth}
      \includegraphics[width=\linewidth]{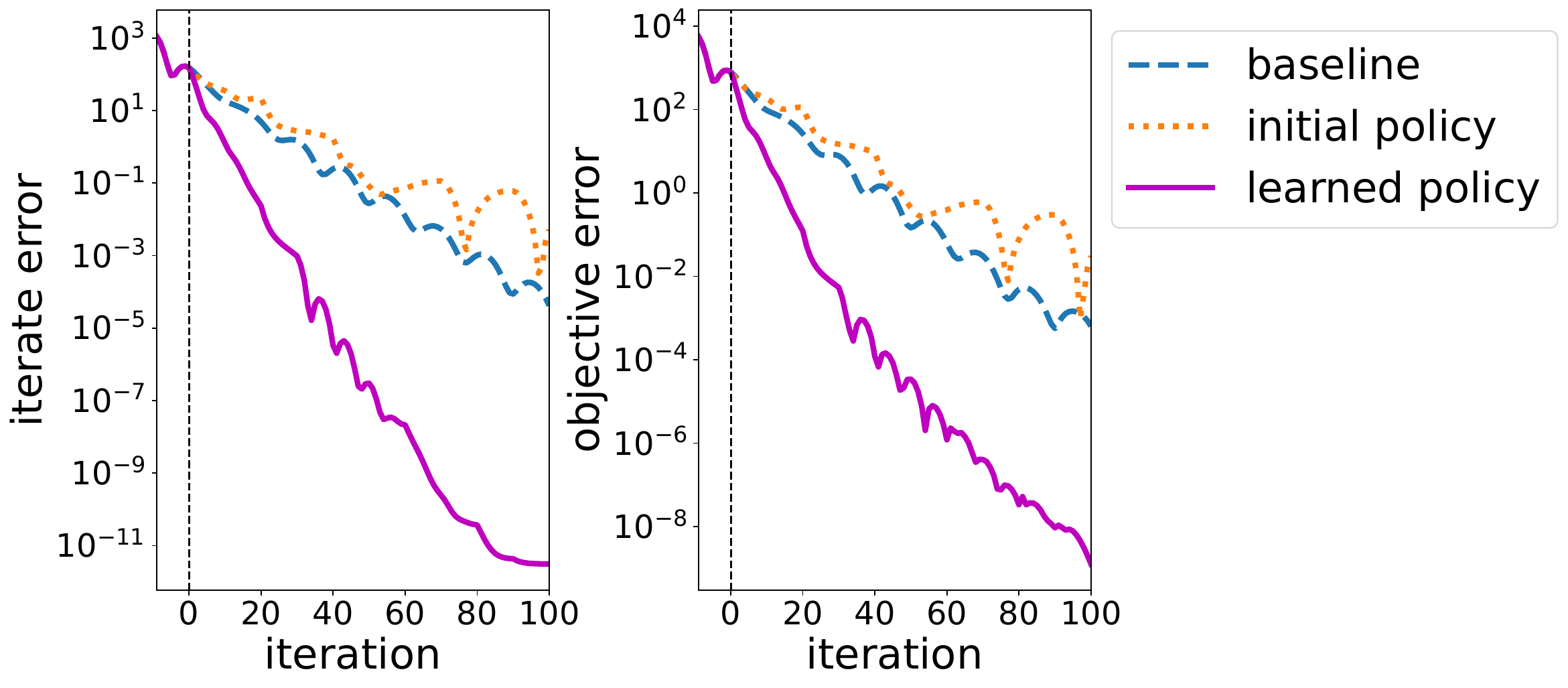}
      \caption{}
      \label{subfig:smooth:train}
    \end{subfigure}
    \begin{subfigure}[b]{0.98\linewidth}
      \includegraphics[width=\linewidth]{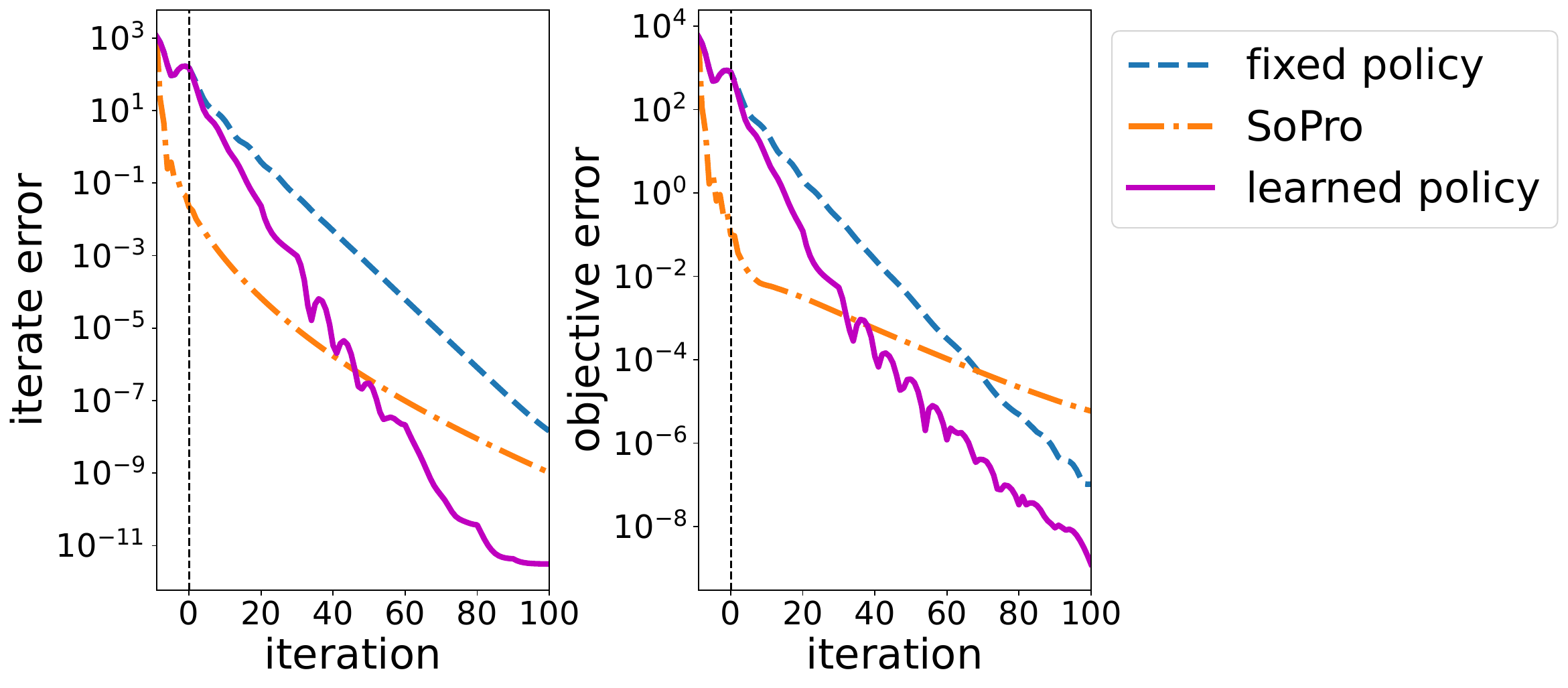}
      \caption{}
      \label{subfig:smooth:compare}
    \end{subfigure}
    \caption{(a) Convergence performance of base model (\ref{eqn:param_DAMM_x})\textminus(\ref{eqn:param_DAMM_q}) under the baseline, the initial policy and the learned policy for solving (\ref{prob:logistic_reg}). (b) Convergence performance of base model (\ref{eqn:param_DAMM_x})\textminus(\ref{eqn:param_DAMM_q}) under the learned policy and the fixed policy (i.e., $\pi(a^c \mid s)\equiv 1$), as well as the convergence performance of SoPro.}
    \label{fig:smooth}
\end{figure}

Figure \ref{subfig:smooth:train} illustrates the training effect by comparing the convergence performance of base model (\ref{eqn:param_DAMM_x})\textminus(\ref{eqn:param_DAMM_q}) under the learned policy with the baseline and the initial policy, where our learning-based framework is able to achieve far more precise solutions within the same number of iteration steps.

We then compare the base model under the learned policy and a hand-selected fixed policy (i.e., $\pi(a^c \mid s)\equiv 1$), respectively, and a recently proposed \textit{second-order} method SoPro \cite{Wu2021SoPro} that is able to address the smooth case of problem (\ref{eqn:min_sum_si_p_ri}) with linear convergence rate and is shown to outperform various existing second-order methods in solving the logistic regression problem \cite{Wu2021SoPro}.
Both the constant action of the fixed policy and SoPro are fine-tuned on the validation set.
The results are shown in Figures \ref{subfig:smooth:compare}.
Still, the learned policy performs better than the fixed policy, indicating the advantage of learning-based adaptive configuration over traditional hand-tuning paradigm.
Since our learning-based framework requires the first $n=10$ iterations to constitute the initial state $s^{(0)}$, the base model (\ref{eqn:param_DAMM_x})\textminus(\ref{eqn:param_DAMM_q}) under the learned policy shows slower descending speed at the beginning compared to SoPro.
However, it catches up quickly as the iterations proceed and finally achieves more accurate solutions.

\subsection{Linear $\ell_1$-Regression with Lasso Regularization} \label{subsec:exp_nonsmooth}

In this subsection, we consider learning a linear model by minimizing the Lasso regularized mean absolute error (MAE) loss over the networked system:
\begin{equation}\label{prob:l1_lasso_reg}
  \min_{x\in\mathbb{R}^d} \sum_{i=1}^{N} \left( \frac{1}{m_i} \sum_{j=1}^{m_i} \lvert a_{i,j}^T x - b_{i,j}\rvert + \lambda \lVert x \rVert_1 \right),
\end{equation}
where $a_{i,j} \in \Rbb^d$ and $b_{i,j}\in\Rbb$ are the feature vector and label of the $j$-th sample on computing node $i$ respectively, $m_i$ is the number of data samples $\{a_{i,j}, b_{i,j}\}$ on each node, and $\lambda > 0$ is the regularization parameter.
This non-smooth objective function is another special case of (\ref{eqn:min_sum_si_p_ri}) with $r_i(x)=\frac{1}{m_i} \sum_{j=1}^{m_i} \lvert a_{i,j}^T x - b_{i,j}\rvert + \lambda \lVert x \rVert_1$ but all $s_i(x)=0$.
Note that in this case both $\nabla s_i(\cdot)$ and $\nabla^2 s_i(\cdot)$ are trivially zero.
Hence, the local observation of each node $i$ in round $t$ simply becomes $o_i^{(t)} = \begin{bmatrix} 
  \sigma_i^{n(t-1)+1}, \cdots, \sigma_i^{nt}
\end{bmatrix}$, and the global action reduces to $a^{(t)} = (\beta^{(t)}, \rho^{(t)})$.
The experiment is also conducted on the real dataset \textit{abalone} as in Section \ref{subsec:exp_composite}.
Similar procedures are adopted to generate the problem instances and form the training set, validation set and test set.

\begin{figure}[t]
  \centering
  \begin{subfigure}[b]{0.99\linewidth}
    \includegraphics[width=\linewidth]{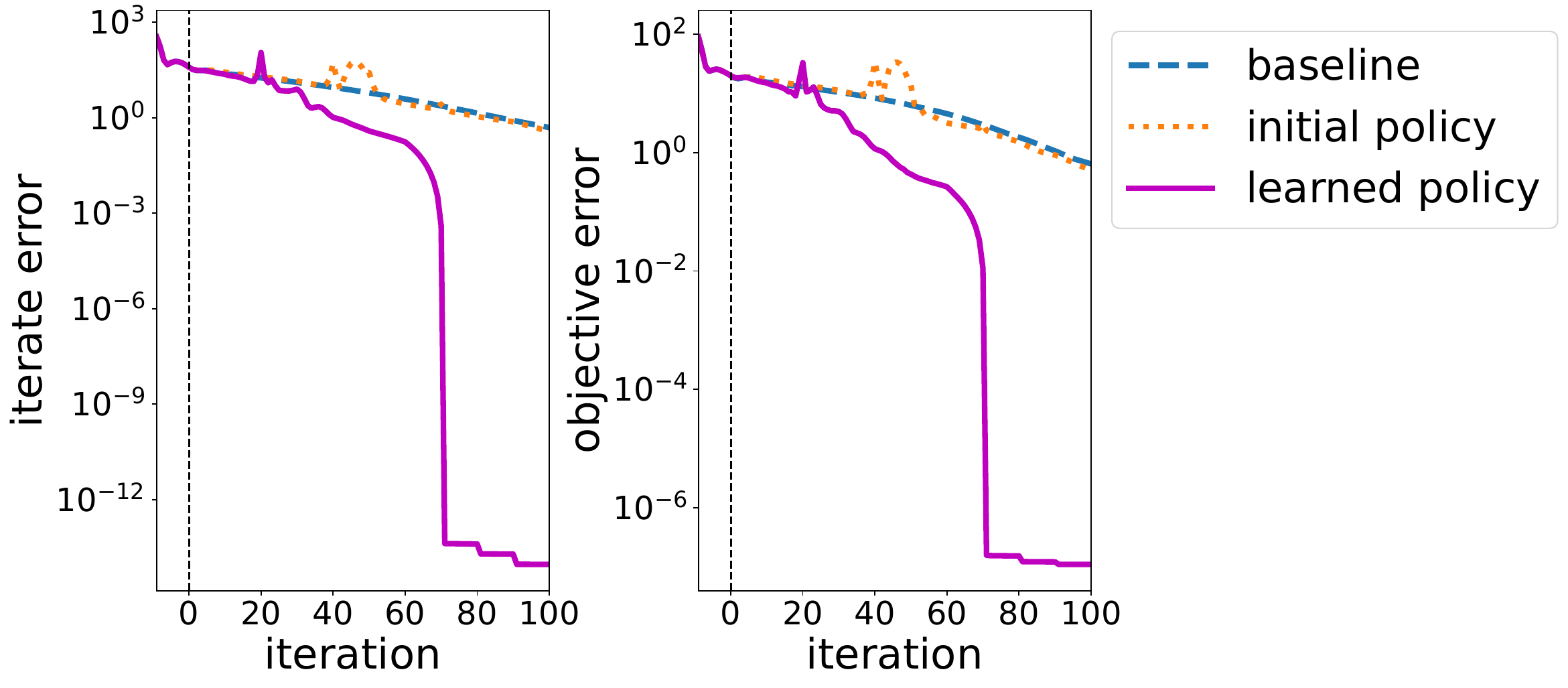}
    \caption{}
    \label{subfig:nonsmooth:train}
  \end{subfigure}
  \begin{subfigure}[b]{0.98\linewidth}
    \includegraphics[width=\linewidth]{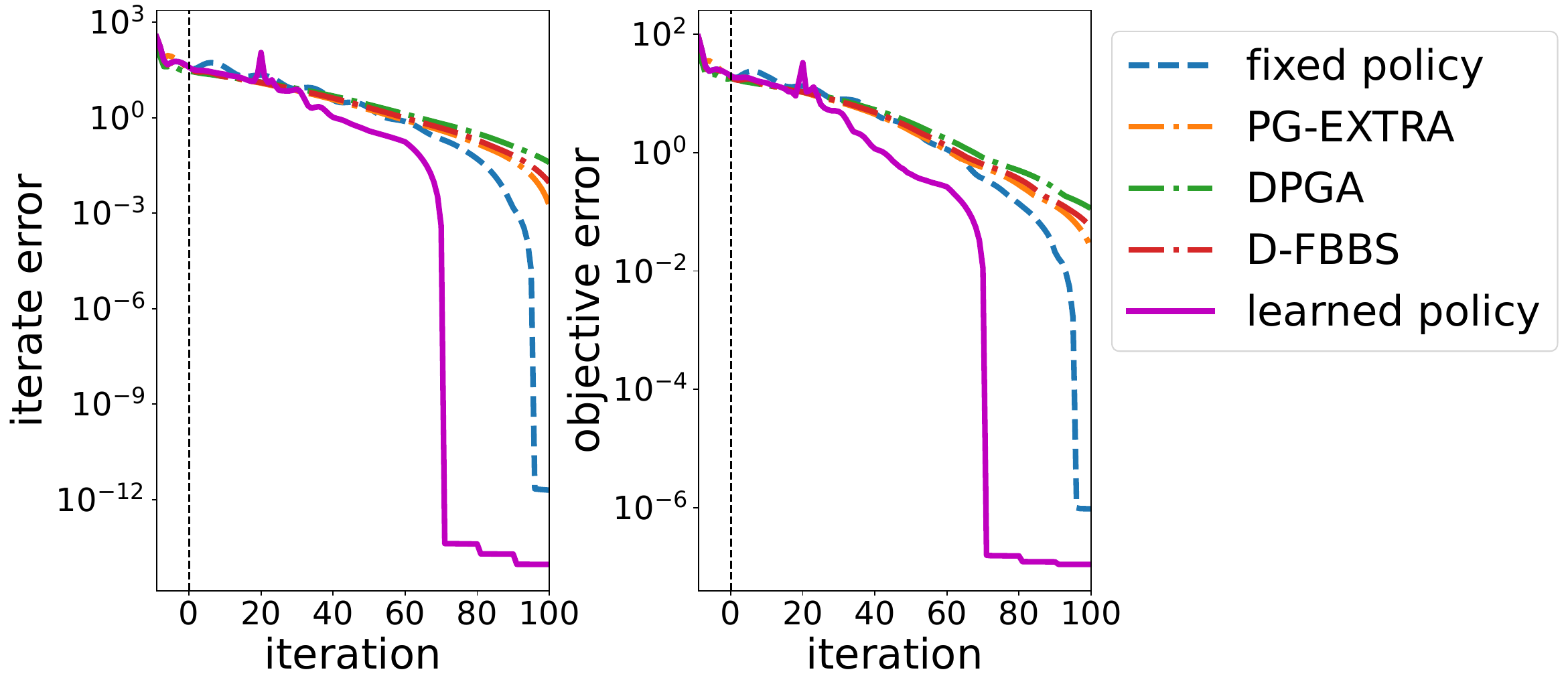}
    \caption{}
    \label{subfig:nonsmooth:compare}
  \end{subfigure}
  \caption{(a) Convergence performance of base model (\ref{eqn:param_DAMM_x})\textminus(\ref{eqn:param_DAMM_q}) under the baseline, the initial policy and the learned policy for solving (\ref{prob:l1_lasso_reg}). (b) Convergence performance of base model (\ref{eqn:param_DAMM_x})\textminus(\ref{eqn:param_DAMM_q}) under the learned policy and the fixed policy (i.e., $\pi(a^c \mid s)\equiv 1$), as well as the convergence performance of state-of-the-art algorithms applicable to (\ref{prob:l1_lasso_reg}).}
  \label{fig:nonsmooth}
\end{figure}

Figure \ref{subfig:nonsmooth:train} compares the convergence performance of the networked system under the baseline, the initial policy and the learned policy in terms of the iterate error and the objective error.
Despite the inferior performance yielded by the baseline or the initial policy, it can be seen that our learning-based framework improves the performance vastly after training.

In Figure \ref{subfig:nonsmooth:compare}, we compare the learned policy with a hand-selected fixed policy (i.e., $\pi(a^c \mid s)\equiv 1$), and, similarly as in Section \ref{subsec:exp_composite}, incorporate PG-EXTRA \cite{Shi2015PGEXTRA}, DPGA \cite{Aybat2018DPGA} and D-FBBS \cite{Xu2018DFBBS} into the comparison to demonstrate the competitive performance of our learning-based framework.
Again, the constant action of the fixed policy and other algorithms for comparison are all fine-tuned using validation instances.
Observe that the learned policy yields faster convergence and more accurate solutions than the fixed policy, and also significantly outperforms PG-EXTRA, DPGA and D-FBBS.
The results are consistent with that in Section \ref{subsec:exp_composite} and \ref{subsec:exp_smooth}, indicating the effectiveness and efficiency of our learning-based method.

\subsection{Generalization to subsequent iterations}

\begin{figure}[t]
  \centering
  \includegraphics[width=0.99\linewidth]{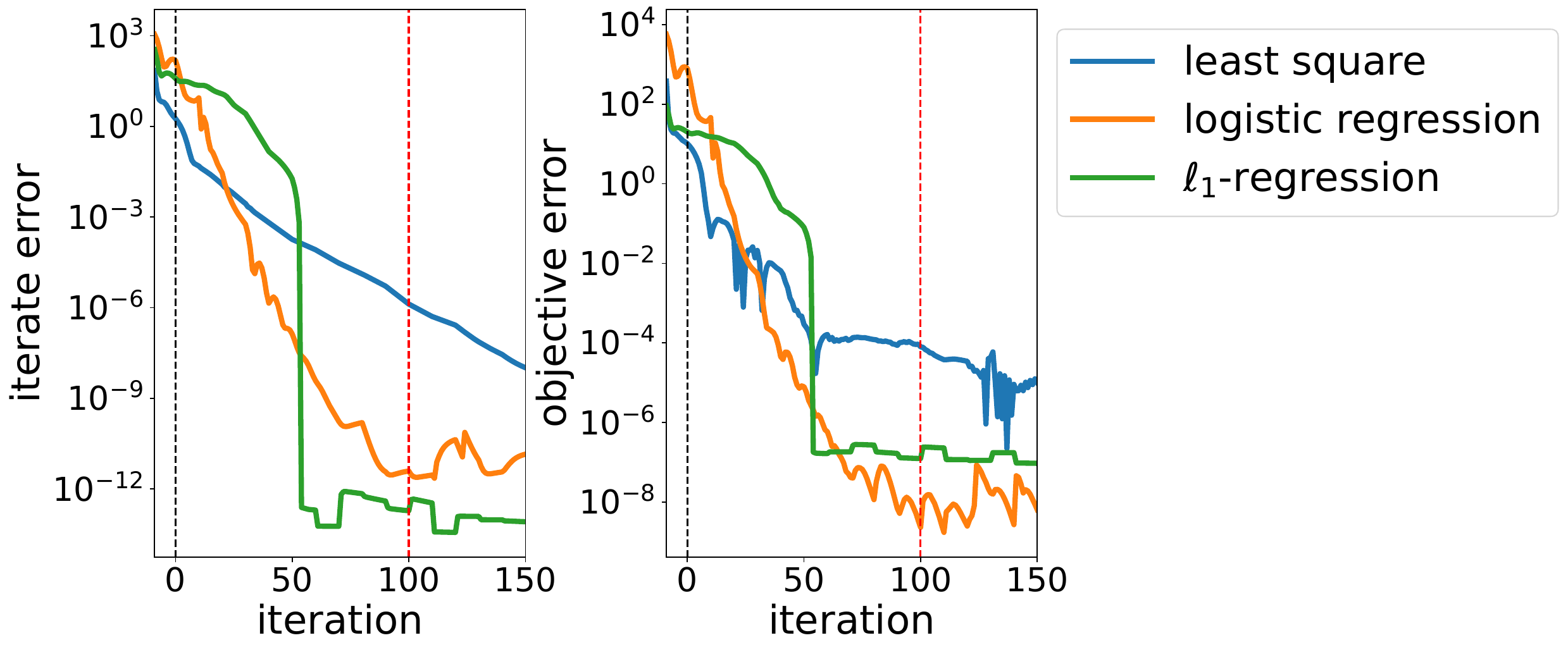}
  \caption{Convergence performance of the networked system under the learned policies for a longer time horizon compared with the training stage. ``Least square" corresponds to the learned policy for solving (\ref{prob:least_square_lasso_reg}), ``logistic regression" corresponds to that for solving (\ref{prob:logistic_reg}), and ``$\ell_1$-regression" corresponds to that for solving (\ref{prob:l1_lasso_reg}). The beginning of the prolonged interval is marked with the vertical dashed line in red.}
  \label{fig:ext_ite}
\end{figure}

This subsection investigates whether our learning-based framework is able to extend its superior performance to subsequent iteration steps, even though it is trained only on a relatively short time horizon.
In practice, this property is important as it provides more flexibility if more iteration steps are needed to further enhance the performance, say, to achieve a more accurate solution.
For each of the problem cases (\ref{prob:least_square_lasso_reg})\textminus(\ref{prob:l1_lasso_reg}), we evaluate the performance of the networked system under the corresponding learned policies for a longer time horizon compared to the training stage.
Figure \ref{fig:ext_ite} shows that all the learned policies for these problem cases are capable of extending their superior performance at least to the subsequent 5 communication rounds, which correspond to $5 \times 10=50$ iteration steps.
For problem case (\ref{prob:least_square_lasso_reg}), the networked system is able to obtain even more accurate solutions with this prolonged time interval.
We note that this observation is essentially disparate with that of \cite{Wang2021}, in which a learning-based optimizer can only iterate as many iteration steps as the number of layers in the constructed GNN, causing its inflexibility in realistic implementation.

\section{Conclusion}
\label{sec:conclusion}
This paper presents a learning-based method for efficient distributed optimization over networked systems.
We have developed a deep reinforcement learning framework that learns a policy for adaptive configuration within a parameterized unifying algorithmic form, which includes various first-order, second-order or hybrid optimization algorithms.
This model-based framework benefits from both the data-driven learning capacity and the convergence guarantee of a theory-driven method.
Adaptive decisions regarding the algorithmic form are made according to local consensus and objective information, which represents the regularities of problem instances and summarizes the solving progress.
We have demonstrated that the proposed method achieves faster convergence and higher accuracy than state-of-the-art distributed optimization algorithms.

\bibliographystyle{IEEEtran}
\bibliography{TCNS-refs}

\begin{IEEEbiography}[]{Daokuan Zhu} received the B.S. degree in electronic information engineering from ShanghaiTech University, Shanghai, China, in 2021.
He is currently pursuing the M.S. degree in information and communication engineering at ShanghaiTech University, Shanghai, China.

His research interests include distributed optimization, federated learning, and learning to optimize.
\end{IEEEbiography}

\begin{IEEEbiography}[]{Tianqi Xu} received the B.S. degree in computer science and technology from ShanghaiTech University, Shanghai, China, in 2022.
She is currently studying toward the M.S. degree in computer science at ShanghaiTech University, Shanghai, China. 

Her research interests include distributed optimization, optimization theory and algorithms, and machine learning. 
\end{IEEEbiography}

\begin{IEEEbiography}[]{Jie Lu} (Member, IEEE) received the B.S. degree in information engineering from Shanghai Jiao Tong University, Shanghai, China, in 2007, and the Ph.D. degree in electrical and computer engineering from the University of Oklahoma, Norman, OK, USA, in 2011. 

She is currently an Associate Professor with the School of Information Science and Technology, ShanghaiTech University, Shanghai, China. Before she joined ShanghaiTech University in 2015, she was a Postdoctoral Researcher with the KTH Royal Institute of Technology, Stockholm, Sweden, and with the Chalmers University of Technology, Gothenburg, Sweden from 2012 to 2015. Her research interests include distributed optimization, optimization theory and algorithms, learning-assisted optimization, and networked dynamical systems.
\end{IEEEbiography}

\end{document}